\documentstyle[amssymb,amsfonts]{amsart}

\def\be#1{\begin{equation} \label{#1}}
\def\bi{\begin{itemize}}
\def\bs{\begin{split}}
\def\es{\end{split}}
\def\ba{\begin{align}}
\def\bas{\begin{align*}}
\def\ea{\end{align}}
\def\eas{\end{align*}}

\def\Re{{\hbox{Re}}}

\def\dist{{\hbox{dist}}}

\def\R{{\hbox{\bf R}}}

\def\T{{\hbox{\bf T}}}

\def\eps{\varepsilon}
\newenvironment{proof}{\noindent {\bf Proof} }{\endprf\par}
\def \endprf{\hfill  {\vrule height6pt width6pt depth0pt}\medskip}
\def\emph#1{{\it #1}}
\def\textbf#1{{\bf #1}}

\theoremstyle{plain}
  \newtheorem{theorem}[subsection]{Theorem}

  \newtheorem{lemma}[subsection]{Lemma}

\theoremstyle{remark}

\theoremstyle{definition}

\include{psfig}

\begin{document}

\title[Orbital instability below the energy norm]{Orbital instability bounds for the non-linear Schr\"odinger equation below the energy norm}
\author{J.~Colliander}
\thanks{J.E.C. is supported in part by N.S.F. grant DMS 0100595 and N.S.E.R.C. grant RGPIN 250233-03.}
\address{University of Toronto }
\author{M.~Keel}
\thanks{M.K. is supported in part by N.S.F. Grant DMS
                         9801558}
\address{University of Minnesota}
\author{G.~Staffilani}
\thanks{G.S. is supported in part by N.S.F. Grant DMS 0100375 and by a grant
from Hewlett and Packard and the Sloan Foundation.}
\address{Stanford University and Brown University}
\author{H.~Takaoka}
\thanks{H.T. is supported in part by J.S.P.S. Grant No. 13740087.}
\address{Hokkaido University}
\author{T.~Tao}
\thanks{T.T. is a Clay Prize Fellow and is supported in part by grants
from the Packard Foundation.}
\address{University of California, Los Angeles}

\vspace{-0.3in}
\begin{abstract}
We continue the study (initiated in \cite{ckstt:7}) of the orbital stability of the ground state cylinder for focusing non-linear Schr\"odinger equations in the $H^s$ norm for $1-\eps < s < 1$.  In the $L^2$-subcritical case we obtain a polynomial bound for the time required to move away from the ground state cylinder.  If one is only in the $H^1$-subcritical case then we cannot show this, but we can obtain global well-posedness and polynomial growth of $H^s$ norms for $s$ sufficiently close to 1. 
\end{abstract}

\maketitle

\section{Introduction}\label{introduction-sec}

We consider the Cauchy problem for the non-linear Schr\"odinger equation
\be{nls}
iu_t + \Delta u = F(u); \quad u(x,0) = u_0(x)
\end{equation}
where $u(x,t)$ is a complex-valued function on $\R^n \times \R$ for some $n \geq 1$, $u_0(x)$ lies in the Sobolev space $H^s(\R^n)$ for some $s \in \R$, and the non-linearity $F(u)$ is the power non-linearity 
$$ F(u) := \pm |u|^{p-1} u$$
for some $p > 1$ and sign $\pm$.  We refer to the + sign as \emph{defocusing} and the - sign as \emph{focusing}.  The long-time behavior of this equation has been extensively studied in the energy class (see e.g. \cite{wein:modulate}, \cite{wein}, \cite{gv:scatter}, \cite{nak:scatter}, \cite{borg:book}), however in this paper we shall be interested in regularities $1 - \eps < s < 1$ slightly weaker than the energy class.  When $p$ is an odd integer, then $F$ is algebraic and there are a number of results addressing the long-time behavior in $H^s$ in this case (\cite{tak:dnls}, \cite{ckstt:1}, \cite{ckstt:4}, \cite{ckstt:5}, \cite{ckstt:6}, \cite{bourg.2d}, \cite{borg:book}, \cite{borg:scatter}).  One of the main purposes of this paper is to demonstrate that these techniques can be partially extended to the non-algebraic case.

We review the known local and global well-posedness theory for this equation.
It is known (see e.g. \cite{cwI}) that the Cauchy problem \eqref{nls} is locally well-posed in $H^s$ when\footnote{When $p$ is not an odd integer one also needs the constraint $\lfloor s \rfloor < p-1$ because of the limited regularity of $F$.  We shall gloss over this technicality since we will primarily be concerned with the regime $0 < s < 1$.} $s \geq \max(0, s_c)$, where $s_c$ is the \emph{critical regularity}
$$ s_c := \frac{n}{2} - \frac{2}{p-1}.$$
In the sub-critical case $s > s_c$ the time of existence depends only on the $H^s$ norm of the initial data.  We refer to the case $s_c < 1$ as the \emph{$H^1$-subcritical case}, and the case $s_c < 0$ as the \emph{$L^2$-subcritical case}. In particular, local well-posedness below $H^1$ is only known in the $H^1$-subcritical case.

These local well-posedness results are known to be sharp in the focusing case (see \cite{sulem}); for instance, one has blowup in arbitrarily small time when $s < s_c$.  In light of the recent work in \cite{cct} it is likely that these results are also sharp in the defocusing case (at least if one wants to have a fairly strong notion of well-posedness). 

In the $H^1$-subcritical case it is known \cite{cwI} that the Cauchy problem \eqref{nls} is globally well-posed in $H^s$ for all $s \geq 1$; indeed, one also has scattering in the defocusing case \cite{gv:scatter}, \cite{nak:scatter}.  The global well-posedness is an easy consequence of the local theory, conservation of the $L^2$ norm and Hamiltonian
\be{hamil-def}
H(u) := \int \frac{1}{2} |\nabla u|^2 \pm \frac{1}{p+1} |u|^{p+1}\ dx,
\end{equation}
combined with the Gagliardo-Nirenberg inequality.
Similarly, in the $L^2$-subcritical case $s_c < 0$ one has global well-posedness in $H^s$ for all $s \geq 0$.

This leaves open the question of the global well-posedness in $H^s$ in the intermediate regime $0 \leq s_c \leq s < 1$.  When the initial data has small norm or is localized in space then global well-posedness is known (see e.g. \cite{sulem}) but the general case is still open\footnote{In the $L^2$-critical case $s_c = 0$ one might expect global well-posedness for large $L^2$ data by combining $L^2$ conservation and the local well-posedness theory.  However a subtlety arises because the $L^2$ norm cannot be scaled to be small, and indeed in the focusing case the $L^2$ mass can concentrate to a point singularity.  In the defocusing case the large $L^2$ global well-posedness remains an important open problem (even in the radial case).}. 

For simplicity we restrict ourselves to the defocusing case. The cases when $p$ is an odd integer have been extensively studied; we summarize the known results in Table \ref{old}.

\begin{figure}\label{old}
\begin{tabular}{|l|l|l|l|} \hline
Dimension $n$ & Exponent $p$ & Critical exponent $s_c$ & Global well-posedness \\
\hline
$1$ & $5$ & $0$ & $s > 32/33$\cite{tak:dnls} \\
&  &  & $s > 2/3$\cite{ckstt:5} \\
 &  &  & $s > 1/2$\cite{ckstt:6} \\
\hline
$2$ & $3$ & $0$ & $s > 2/3$\cite{bourg.2d}, \cite{borg:book} \\
 &  &  & $s > 3/5$\cite{bourg.2d} \\
 &  &  & $s > 4/7$\cite{ckstt:1} \\
 &  &  & $s > 1/2$\cite{ckstt:4} \\
\hline
$3$ & $3$ & $1/2$ & $s > 11/13$\cite{borg:scatter} \\
 &  &  & $s > 5/6$\cite{ckstt:1} \\
 &  &  & $s > 5/7$ (radial data only) \cite{borg:scatter} \\
&  &  & $s > 4/5$\cite{ckstt:scatter} \\
\hline
\end{tabular}
\caption{Known global well-posedness results in the regime $0 \leq s_c \leq s < 1$. In these results the non-linearity is assumed to be defocusing.}
\end{figure}

Our first main result to extend these results\footnote{Similar results have been obtained for the non-linear wave equation in \cite{kpv:nlw}.  However the argument in \cite{kpv:nlw} cannot be extended to NLS because it relies on the gain in regularity inherent in the wave Strichartz estimates, which are not present for the Schr\"odinger Strichartz estimates.}
 to fractional $p$ in the regime $0 \leq s_c < 1$.

\begin{theorem}\label{main-1}
Suppose that we are in the defocusing case with $s_c < 1$.  Then the Cauchy problem \eqref{nls} is globally well-posed in $H^s$ whenever
$$
1 > s > 1 - \eps(n,p)
$$
for some $\eps(n,p) > 0$.
Furthermore, one has the polynomial growth bound
\be{poly}
\| u(T) \|_{H^s} \leq C(s, n, p, \|u_0\|_{H^s}) (1 + |T|)^{C(s,n,p)}
\end{equation}
for all times $T \in \R$.
\end{theorem}

One can probably modify these arguments and exponents to handle the focusing case, especially in the $L^2$-subcritical case $s_c \leq 0$, but we shall not do so here.  In the $L^2$-subcritical case one has global well-posedness for all $s \geq 0$ thanks to $L^2$ norm conservation, but the question of polynomial growth of $H^s$ norms for $s$ close to 0 remains open. 

We have not explicitly calculated $\eps(n,p)$; the exponents given by our arguments are significantly weaker than those in the results previously mentioned.  

Our approach is based on the ``$I$-method'' in \cite{keel:mkg}, \cite{ckstt:1}, \cite{ckstt:2}, \cite{ckstt:3}, \cite{ckstt:4}, \cite{ckstt:5}, \cite{ckstt:6} (see also \cite{keel:wavemap}).  The idea is to replace the conserved quantity $H(u)$, which is no longer available when $s < 1$, with an ``almost conserved'' variant $H(Iu)$, where $I$ is a smoothing operator of order $1-s$ which behaves like the identity for low frequencies (the exact definition of ``low frequencies'' will depend ultimately on the time $T$). Since $p$ is not necessarily an odd integer, we cannot use the multi-linear calculus (or $X^{s,b}$ spaces) in previous papers, and must rely instead on more rudimentary tools such as Taylor expansion (and Strichartz spaces).  In particular there does not appear to be an easy way to improve the exponent $\eps(n,p)$ by adding correction terms to $H(Iu)$ (cf. \cite{ckstt:2}, \cite{ckstt:4}, \cite{ckstt:6}).  Also, we will avoid the use of $L^2$ conservation law as much as possible as this norm can be critical or supercritical.  As a partial substitute we shall use the subcritical $L^{p+1}$ norm which we can control from \eqref{hamil-def} (cf. \cite{keel:mkg}).  

In the cases when $p$ is an odd power, smoothing estimates such as the bilinear Strichartz estimate of Bourgain (see e.g. \cite{borg:book}) are very useful for these types of results.  However we will not use any sort of smoothing estimates in our analysis\footnote{Indeed, the sharp bilinear improvements to Strichartz' inequality in $L^p_{x,t}$ norms are only known in dimensions $n=1,2$; for higher dimensions the problem is connected with the very difficult restriction problem for the paraboloid (see e.g. \cite{tvv:bilinear} for a discussion).}, and rely purely on Strichartz estimates instead.  (One of the advantages of the $I$-method is that one can use commutator estimates involving the operator $I$ as a substitute for smoothing estimates even when the nonlinearity has no smoothing properties).

Our second result concerns the orbital stability of ground states in the focusing case.  For this result we shall restrict ourselves\footnote{In the $L^2$-critical or $L^2$-supercritical cases the ground state is known to be unstable, indeed one can have blowup in finite time even for data arbitrarily close to a ground state in smooth norms.  See \cite{merle}.} to the $L^2$-subcritical case $s_c < 0$. 

When $s_c < 0$ there exists a unique radial positive Schwartz function $Q(x)$ which solves the equation
$\Delta Q - Q = F(Q)$ 
on $\R^n$ (see \cite{coff}, \cite{mc}, \cite{kwong}).  We refer to $Q$ as the \emph{canonical ground state} at energy\footnote{Other energies $E$ are possible but can be easily obtained from the energy 1 state by scaling.} 1.  The Cauchy problem \eqref{nls} with initial data $u_0 = Q$ then has an explicit solution $u(t) = e^{it} Q$.  More generally, for any $x_0 \in \R^n$ and $e^{i\theta} \in S^1$, the Cauchy problem with initial data $u_0(x) = e^{i\theta} Q(x-x_0)$ has explicit solution $e^{i(\theta+t)} Q(x-x_0)$.  If we thus define the \emph{ground state cylinder} $\Sigma$ by
$$ \Sigma := \{ e^{i\theta} Q(\cdot - x_0): x_0 \in \R^n, e^{i\theta} \}$$
we see that the non-linear flow \eqref{nls} preserves $\Sigma$.  We now investigate how the non-linear flow \eqref{nls} behaves on neighborhoods of $\Sigma$.

In \cite{wein:modulate}, \cite{wein} Weinstein showed that in the $L^2$-subcritical case $s_c < 0$ and when $n=1,3$, the ground state cylinder $\Sigma$ was $H^1$-stable.  More precisely, he showed an estimate of the form
\be{stab}
\dist_{H^1}(u(t), \Sigma) \sim \dist_{H^1}(u_0, \Sigma)
\end{equation}
for all solutions $u_0$ to \eqref{nls} and all times $t \in \R$.  In other words, solutions which started close a ground state in $H^1$ stayed close to a ground state for all time (though the nearby ground state may vary in time).

To prove \eqref{stab}, Weinstein employed the \emph{Lyapunov functional}
\be{L-def}
L(u) := 2H(u) + \int |u|^2
= \int |\nabla u|^2 + |u|^2 + \frac{2}{p+1} |u|^{p+1}\ dx,
\end{equation}
which is well-defined for all $u \in H^1$.  Since this quantity is a combination of the Hamiltonian \eqref{hamil-def} and the $L^2$ norm, it is clearly an invariant of the flow \eqref{nls}.  The ground states in $\Sigma$ then turn out to minimize $L$, so that $L(u) \geq L(Q)$ for all $u \in H^1$.  More precisely, we have the inequality
\be{weinstein}
L(u) - L(Q) \sim \dist_{H^1}(u(t), \Sigma)^2 \hbox{ whenever } \dist_{H^1}(u(t), \Sigma) \ll 1;
\end{equation}
see \cite{wein:modulate}, \cite{wein}.  The stability estimate \eqref{stab} then follows easily from \eqref{weinstein} and the conservation of $L$.

Weinstein's proof of \eqref{weinstein} requires the uniqueness of the ground state $Q$, which at the time was only proven for $n=1,3$ \cite{coff}.  However, this uniqueness result has since been extended to all dimensions $n$ \cite{kwong} (with an earlier partial result in \cite{mc}).  Thus \eqref{weinstein} (and hence \eqref{stab}) holds for all dimensions $n$ (always assuming that we are in the $L^2$ subcritical case $s_c < 0$, of course).

In \cite{ckstt:7} the $H^1$ orbital stability result was partially extended to regularities $H^s$, $0 < s < 1$, in the special case $n=1$, $p=3$ (which among other things is completely integrable).  The second main result of this paper is a partial extension of the results of \cite{ckstt:7} to arbitrary $L^2$-subcritical NLS: 

\begin{theorem}\label{main-2}  Consider the focusing case with $s_c < 0$.  Suppose that $u(t)$ solves \eqref{nls} with
$\dist_{H^s}(u_0, \Sigma) \lesssim \sigma$ for some $0 < \sigma \ll 1$ and $1 - \eps(n,p) < s < 1$.  Then one has
$\| u(t) \|_{H^s} \lesssim 1$ for all $t = O(\sigma^{-C(n,p,s)})$ for some $C(n,p,s) > 0$.
\end{theorem}

In other words, if the initial data stays close to the ground state cylinder in $H^s$ norm, then the solution stays inside a ball of bounded radius in $H^s$ for a fairly long period of time.  (After this time, one can use Theorem \ref{main-1} to give polynomial bounds on the growth of the $H^s$ norm).  Note that if one were to try to naively use perturbation theory to prove this theorem, one would only be able to keep $u(t)$ inside this ball for times $t = O(\log(1/\sigma))$ (because after each time interval of length $\sim 1$, the distance to the ground state cylinder might conceivably double).

Of course, when $s = 0$ or $s = 1$ one can use the conservation laws to obtain Theorem \ref{main-2} for all time $t$.  However there does not seem to be any easy way to interpolate these endpoint results to cover the $0 < s < 1$ case, since the flow \eqref{nls} is neither linear nor complex-analytic. 

The proof of Theorem \ref{main-2} also proceeds via the ``I-method''.  The main idea is to show that the modified Lyapunov functional $L(Iu)$ is ``almost conserved''.  It should be possible to refine this method by approximating $u(t)$ carefully by a ground state as in \cite{ckstt:7} and obtain a more precise estimate of the form
$$ \dist(u(t),\Sigma) \lesssim (1+|t|)^{C(n,p,s)} \dist(u_0,\Sigma)$$
for all time $t$.  (In \cite{ckstt:7} this is achieved in the model case $n=1$, $p=3$).  However there seem to be some technical difficulties in making this approach viable, mainly due to the lack of regularity\footnote{The specific obstacle is as follows.  In our current argument we must estimate commutator expressions such as $IF(u) - F(Iu)$.  To utilize the ground state cylinder as in \cite{ckstt:7} one would also consider expressions such as $F(Q+w)-F(Q)$.  To use both simultaneously one needs to estimate a double difference such as $I(F(Q+w)-F(Q)) - (F(Q+Iw)-F(Q))$.  However when $p<2$, $F$ is not twice differentiable, and so correct estimation of the double difference seems very subtle.} of the non-linearity $F$, and so we will not pursue this matter. 

It should be possible to remove the constraint $s > 1 - \eps(n,p)$ and prove Theorem \ref{main-2} for all $s > 0$ (as in \cite{ckstt:7}).  This may however require some additional assumptions on $p$ (e.g. one may need $p < 1 + \frac{2}{n}$) as it becomes difficult to control the modified Hamiltonian for $s$ close to zero otherwise. 

The authors thank Monica Visan for some helpful corrections.

\section{Preliminaries: Notation}\label{notation-sec}

Throughout the paper $n$, $p$, $s$ are considered to be fixed.  We will always have the implicit assumption ``$1 - \eps(n,p) < s < 1$ for some sufficiently small $\eps(n,p) > 0$'' in our arguments.
We let $A \lesssim B$ or $A = O(B)$ denote the estimate $A \leq C B$, where $C$ is a positive constant which depends only on $n$, $p$, $s$.

We compute the derivatives $F_z$, $F_{\overline z}$ of $F(z) = \pm |z|^{p-1} z$
as
$$ F_z(z) := \pm \frac{p+1}{2} |z|^{p-1}; \quad F_{\overline{z}}(z) := \pm \frac{p-1}{2} |z|^{p-3} z^2.$$
We write $F'$ for the vector $(F_z, F_{\overline{z}})$, and adopt the notation
$$ w \cdot F'(z) := w F_z(z) + \overline{w} F_{\overline{z}}(z).$$
In particular we observe the chain rule
$$ \nabla F(u) = \nabla u \cdot F'(u).$$

Clearly $F'(z) = O(|z|^{p-1})$. We also observe the useful H\"older continuity estimate
\be{holder}
|F'(z) - F'(w)| \lesssim |z-w|^\theta (|z|^{p-1-\theta} + |w|^{p-1-\theta})
\end{equation}
for all complex $z,w$, where $\theta := \min(p-1,1)$ is a number in the interval $(0,1]$.  In a similar spirit we record the estimate 
\be{triangle}
F(z+w) = F(z) + O(|w| |z|^{p-1}) + O(|w|^p)
\end{equation}
for all complex $z,w$.

We define the spatial Fourier transform on $\R^n$ by
$$ \hat f(\xi) := \int e^{-2\pi i x \cdot \xi} f(x)\ dx.$$

For $N > 1$, we define the Fourier multiplier $I = I_N$ by
$$ \widehat{Iu}(\xi) := m(\xi/N) \hat u(\xi)$$
where $m$ is a smooth radial function which equals 1 for $|\xi| \leq 1$ and equals $|\xi|^{s-1}$ for $|\xi| \geq 2$.  Thus $I_N$ is an operator which behaves like the {\bf I}dentity for low frequencies $|\xi| \leq N$, and behaves like a (normalized) {\bf I}ntegration operator of order $1-s$ for high frequencies $|\xi| \gtrsim N$.  In particular, $I$ maps $H^s$ to $H^1$ (but with a large operator norm, roughly $N^{1-s}$).  This operator will be crucial in allowing us to access the $H^1$ theory at the regularity of $H^s$.  We make the useful observation that $I$ has a bounded convolution kernel and is therefore bounded on every translation-invariant Banach space.

We also define the fractional differentiation operators $|\nabla|^\alpha$ for real $\alpha$ by
$$ \widehat{|\nabla|^\alpha u}(\xi) := |\xi|^\alpha \hat u(\xi)$$
and the modified fractional differentiation operators $\langle \nabla \rangle^\alpha$ by
$$ \widehat{\langle \nabla\rangle^\alpha u}(\xi) := \langle \xi\rangle^{\alpha} \hat u(\xi)$$
where $\langle \xi \rangle := (1 + |\xi|^2)^{1/2}$.

We then define the inhomogeneous Sobolev spaces $H^s$ and homogeneous Sobolev spaces $\dot H^s$ by
$$ \| u \|_{H^s} := \| \langle \nabla \rangle^s u \|_{L^2_x}; \quad \| u \|_{\dot H^s} := \| |\nabla|^s u \|_{L^2_x}.$$

We shall frequently use the fact (from elementary Littlewood-Paley  theory, see e.g. \cite{stein:small}) that if $u$ has Fourier transform supported on a set $|\xi| \lesssim M$, then one can freely replace positive powers of derivatives $\nabla$ by the corresponding powers of $M$ in $L^p$ norms for $1 < p < \infty$, thus for instance
$$ \| \langle \nabla \rangle^s u \|_p \lesssim \langle M \rangle^s \| u \|_p$$
for $s > 0$.  Conversely, if $u$ has Fourier transform supported on $|\xi| \gtrsim M$ then one can insert positive powers derivatives $\nabla$ and gain a negative power of $M$ as a result:
$$ \| u \|_p \lesssim \langle M \rangle^{-s} \| \langle \nabla \rangle^s u \|_p.$$
In particular, for any $u$, $u-Iu$ has Fourier support on the region $|\xi| \geq N$, hence
$$ \| u - Iu \|_p \lesssim N^{-\eps} \| \langle \nabla \rangle^\eps u \|_p$$
for any $\eps > 0$.  This fact (and others similar to it) will be key in extracting crucial negative powers of $N$ in our estimates.

\section{Preliminaries: Strichartz spaces}\label{strichartz-sec}

In this section we introduce the $H^1$ Strichartz spaces we will use for the semilinear equation \eqref{nls}, and derive the necessary nonlinear estimates for our analysis.  In particular, we obtain nonlinear commutator estimates involving the fractional nonlinearity $F(u)$, which is the main new technical advance in this paper.

We will always assume we are in the $H^1$-subcritical case
\be{h1-subcrit}
\frac{1}{p-1} > \frac{n-2}{4}.
\end{equation}

Let $t_0$ be a time and $0 < \delta \leq 1$.  In what follows we restrict spacetime to the slab $\R^n \times [t_0, t_0 + \delta]$.  
We define the spacetime norms $L^q_t L^r_x$ by
$$ \| u \|_{L^q_t L^r_x} := (\int \| u(t) \|_{L^r_x}^q\ dt)^{1/q}.$$
We shall often abbreviate $\| u\|_{L^q_t L^r_x}$ as $\| u \|_{q,r}$.

We shall need a space $L^{q_0}_t L^{r_0}_x$ to hold the solution $u$, another space $L^{q_1}_t L^{r_1}_x$ to hold the derivative $\nabla u$ (or $I \nabla u$), and the dual space $L^{q'_1}_t L^{r'_1}_x$ to hold the derivative non-linearity $\nabla F(u)$ (or $I\nabla F(u)$).  (Here of course $q'$ denotes the exponent such that $1/q + 1/q' = 1$).  We also need a space $L^{q_0/(p-1)}_t L^{r_1/(p-1)}_x$ to hold $F'(u)$.  To choose the four exponents $q_0, r_0, q_1, r_1$ we use the following lemma (cf. \cite{cwI}):

\begin{lemma}\label{exponents}  Let $p$ be as above.  There exist exponents $2 < q_0, r_0, q_0, r_1 < \infty$ and $0 < \beta < 1$ obeying the relations
\begin{align}
\frac{2}{q_1} + \frac{n}{r_1} &= \frac{n}{2} \label{x-scaling}\\
\frac{2}{q_0} + \frac{n}{r_0} &= \frac{n-2}{2} + \beta \label{x-scaling-2}\\
\frac{1}{r_1} + \frac{p-1}{r_0} &= \frac{1}{r'_1} \label{x-gap}\\
\frac{1}{q_1} + \frac{p-1}{q_0} &< \frac{1}{q'_1} \label{x-gap-2}\\
r_0 &> p+1.\label{x-hamil}
\end{align}
\end{lemma}

\begin{proof}
We first choose $\beta$ such that so that
$$ \max(0, \frac{n-2}{2}) < \frac{n-2}{2}+\beta < \min(\frac{n}{2}, \frac{2}{p-1});$$
such an $\beta$ exists from \eqref{h1-subcrit}.  Next, we choose $2 < q_0 < \infty$ and $p+1 <  r_0 < \infty$ so that \eqref{x-scaling-2} holds; 
such a pair $q_0, r_0$ exists since
$$ \frac{2}{2} + \frac{n}{p+1} \geq \frac{n}{2} > \frac{n-2}{2}+\beta > 0.$$
Next, we choose $r_1$ so that \eqref{x-gap} holds.  Finally we choose $q_1$ so that \eqref{x-scaling} holds.

From construction we see that \eqref{x-scaling}, \eqref{x-scaling-2}, \eqref{x-gap}, \eqref{x-hamil} hold, and that $2 < q_0, r_0 < \infty$ and $0 < \beta < 1$.  To prove \eqref{x-gap-2}, we observe from \eqref{x-scaling} that
$$ \frac{2}{q'_1} + \frac{n}{r'_1} = \frac{n+4}{2}$$
so by \eqref{x-scaling}, \eqref{x-gap}, and \eqref{x-scaling-2} it suffices to show that
$$ \frac{n}{2} + (p-1)(\frac{n-2}{2}+\beta) < \frac{n+4}{2}.$$
But this follows since $\frac{n-2}{2}+\beta < \frac{2}{p-1}$.

Since $p+1 < r_0 < \infty$, we see from \eqref{x-gap} that $2 < r_1 < p+1$.  In particular we have
$$ \frac{1}{2} > \frac{1}{r_1} > \frac{1}{p+1} > \frac{n-2}{2n}$$
so by \eqref{x-scaling} we have $2 < q_1 < \infty$.  All the required properties are thus obeyed.
\end{proof}

Henceforth $q_0, r_0, q_1, r_1$ are assumed to be chosen as above. 
We define the spacetime norm $X$ by
\be{X-def}
\| u \|_X := \| \nabla u(t_0) \|_2 + \| \nabla (iu_t + \Delta u) \|_{q'_1, r'_1}
\end{equation}

This $X$ norm looks rather artificial, but it is easy to estimate for solutions of \eqref{nls}, and it can be used to control various spacetime norms of $u$.  Indeed, we recall that the hypotheses \eqref{x-scaling}, \eqref{x-scaling-2} imply the (scale-invariant) Strichartz estimate (see e.g. \cite{tao:keel})

\begin{lemma}\label{strich}  We have
$$ \| \nabla u \|_{q_1, r_1} + \| |\nabla|^\beta u \|_{q_0, r_0} \lesssim \| u \|_X.$$
\end{lemma}

In future applications we shall need to control $\| u \|_{q_0, r_0}$ and $\| u \|_{2p, 2p}$ in addition to the norms already controlled by Lemma \ref{strich}.  These norms cannot be controlled purely by the $X$ norm, however from the conservation of the Hamiltonian we will also be able to control\footnote{We also have $L^2$ norm conservation which gives bounds on $L^\infty_t L^2_x$, but it turns out that in the $L^2$-supercritical case these bounds are not favorable.} the $L^\infty_t L^{p+1}_x$ norm, and by combining these estimates we shall be able to estimate everything we need.  More precisely, we have

\begin{lemma}\label{control}  Suppose that we are working on the slab $\R^n \times [t_0, t_0+\delta]$ and that $u$ is a function on this slab obeying the estimates
$$
\| u \|_X + \| u \|_{\infty, p+1} \lesssim 1.
$$
Then we have
\be{beta}
\| \langle \nabla \rangle^\beta u \|_{q_0, r_0} \lesssim \| u \|_X^\kappa
\end{equation}
\be{gamma}
\| \langle \nabla \rangle^\gamma u\|_{2p, 2p} \lesssim \| u \|_X^\kappa
\end{equation}
for some $\gamma, \kappa > 0$.
\end{lemma}

\begin{proof}  We introduce the frequency cutoff $\lambda := \| u \|_X^\eps$ for some $\eps > 0$ to be determined later, and smoothly divide $u = u_{low} + u_{high}$, where $u_{low}$ has Fourier support in the region $|\xi| \lesssim \lambda$ and $u_{high}$ has Fourier support in the region $|\xi| \gtrsim \lambda$.

Consider the contribution of $u_{high}$.  Then $\langle \nabla \rangle$ is bounded by $\lambda^{-1} |\nabla|$, and so \eqref{beta} follows from Lemma \ref{strich}.  To prove \eqref{gamma}, we let $r$ be the exponent such that
$$ \frac{2}{2p} + \frac{n}{r} = \frac{n}{2}.$$
Observe that $2 < r < \infty$.  From Strichartz' estimate (see e.g. \cite{tao:keel}) we have
$$ \| \nabla u \|_{2p,r} \lesssim \| u \|_X.$$
From \eqref{h1-subcrit} we have
$$ \frac{n}{r} - \frac{n}{2p} = \frac{n}{2} - \frac{n+2}{2p} < 1.$$
The claim \eqref{gamma} then follows from Sobolev embedding and the high frequency assumption $|\xi| \gtrsim \lambda$.

Now consider the contribution of $u_{low}$.  Then the $\langle \nabla \rangle$ can be discarded.  Since $r_0$ and $2p$ are both strictly larger than $p+1$, we see from Bernstein's inequality (or Sobolev embedding) that
$$ \| u_{low} \|_{\infty, r_0}, \| u_{low} \|_{\infty, 2p} \lesssim \lambda^c \| u_{low} \|_{\infty, p+1}$$
for some $c>0$.  But by hypothesis we have $\| u_{low} \|_{\infty, p+1} \lesssim 1$.  The claim then follows from a H\"older in time.
\end{proof}

We now use these estimates to prove some non-linear estimates involving $F$ and $I$.  We begin with a bound on $F'(Iu)$.

\begin{lemma}\label{frac-chain}  Suppose that
$$
\| Iu \|_X + \| Iu \|_{\infty, p+1} \lesssim 1.
$$
Then
$$
\| \langle \nabla \rangle^\eps F'(u) \|_{q_0/(p-1), r_0/(p-1)} \lesssim \| u \|_X^{(p-1)\kappa}
$$
for some $0 < \eps = \eps(p,n) < 1$.
\end{lemma}

\begin{proof}
From Lemma \ref{control} we have
\be{u-c}
\| \langle \nabla \rangle^{\beta-\eps} u \|_{q_0, r_0} \lesssim
\| \langle \nabla \rangle^\beta I^{-1} u \|_{q_0, r_0} \lesssim \| u \|_X^{\kappa}.
\end{equation}
In particular we have
$$ \| F'(u) \|_{q_0/(p-1), r_0/(p-1)} \lesssim 
\| u \|_{q_0, r_0}^{p-1} \lesssim \| u \|_X^{(p-1)\kappa}.$$
To get the additional $\eps$ of regularity we shall use H\"older norms.
Since Sobolev norms control H\"older norms (see e.g. \cite{stein:small}) we have have 
$$ \| u - u_y \|_{q_0, r_0} \lesssim |y|^\beta \| u \|_X^{\kappa}
$$
for all $|y| \lesssim 1$, where $u_y(x,t) := u(x-y,t)$ is the translation of $u$ in space by $y$.  From \eqref{holder} and H\"older's inequality we thus have
$$ \| F'(u) - F'(u_y) \|_{q_0/(p-1), r_0/(p-1)}
\lesssim \| u - u_y \|_{q_0, r_0}^\theta (
\| u \|_{q_0, r_0}^{p-1-\theta} + 
\| u_y \|_{q_0, r_0}^{p-1-\theta})
$$
where $\theta := \min(p-1,1)$.
Using \eqref{u-c} and the observation that $F'(u_y) = F'(u)_y$, we thus have a H\"older bound on $F'(u)$:
$$ \| F'(u) - F'(u)_y \|_{q_0/(p-1), r_0/(p-1)} \lesssim |y|^{\beta \theta} \| u \|_X^{(p-1)\kappa}.$$
This yields the desired Sobolev regularity bound for any $0 < \eps < \beta \theta$ (see \cite{stein:small}).
\end{proof}

From this Lemma one can already recover the proof (from \cite{cwI}) of $\dot H^1$ local well-posedness of the NLS equation \eqref{nls}.  Indeed, if $u$ solves \eqref{nls}, then from \eqref{X-def} and the chain rule we have
$$ \| u \|_X \lesssim \| u(t_0) \|_{\dot H^1} + \| \nabla u \cdot F'(u) \|_{q'_1,r'_1},$$
which by H\"older, \eqref{x-gap}, \eqref{x-gap-2}, and Lemma \ref{frac-chain} (discarding the epsilon gain of regularity) yields
$$ \| u \|_X \lesssim \| u(t_0) \|_{\dot H^1} + \| u \|_X^{1+(p-1)\kappa},$$
which (together with a similar inequality for differences in iterates of \eqref{nls}) allows one to obtain well-posedness if the $\dot H^1$ norm of the initial data is sufficiently small.  (Large data can then be handled by a scaling argument).  

By a variant of the argument just described, we can obtain bounds on $I\nabla F(u)$ and the related commutator expression $\nabla (IF(u) - F(Iu))$:

\begin{lemma}\label{chain-rule} Suppose that we are working on the slab $\R^n \times [t_0, t_0+\delta]$ and that $u$ is a function on this slab obeying the estimates
\be{small-x}
\| Iu \|_X + \| Iu \|_{\infty, p+1} \lesssim 1.
\end{equation}
  Then
\be{chain}
\| I \nabla F(u) \|_{q'_1, r'_1} \lesssim \delta^c \| Iu \|_X^{1 + (p-1)\kappa}
\end{equation}
for some $c > 0$.  Furthermore, we have
\be{chain-2}
\| \nabla (I F(u) - F(Iu)) \|_{q'_1, r'_1} \lesssim N^{-\alpha}
\delta^c \| Iu \|_X^{1 + (p-1)\kappa}
\end{equation}
for some $\alpha > 0$.  (The quantities $c, \alpha$ depend of course on $n$, $p$, $s$, $\eps$).
\end{lemma}

One can think of \eqref{chain} as a type of fractional chain rule for the differentiation operator $I \nabla$.
The additional gain of $N^{-\alpha}$ in \eqref{chain-2} arises from the spare epsilon of regularity in Lemma \ref{frac-chain} and the fact that $I$ is the identity for frequencies $\lesssim N$; this gain is crucial to all the results in this paper. 

\begin{proof}
From Lemma \ref{strich} and Lemma \ref{frac-chain} we have
\be{a}
\| Iu \|_{q_1, r_1} \lesssim \| Iu \|_X
\end{equation}
and
\be{b}
\| \langle \nabla \rangle^\eps F'(u) \|_{q_0/(p-1), r_0/(p-1)} \lesssim \| Iu \|_X^{(p-1)\kappa}.
\end{equation}
To utilize \eqref{a}, \eqref{b} we use the following bilinear estimates.

\begin{lemma}\label{comm}  If $s$ is sufficiently close to 1 (depending on $\eps$), then we have
$$ \| I(fg) \|_{q'_1, r'_1} \lesssim \delta^c 
\| I f \|_{q_1, r_1} 
\| \langle \nabla \rangle^\eps g \|_{q_0/(p-1), r_0/(p-1)}$$
and
$$ \| I(fg) - (If) g \|_{q'_1, r'_1} \lesssim \delta^c N^{-\alpha} \| I f \|_{q_1, r_1} \| \langle \nabla \rangle^\eps g \|_{q_0/(p-1), r_0/(p-1)}$$
for some $\alpha > 0$ depending on $s$ and $\eps$, and any $f, g$ on the slab $\R^n \times [t_0, t_0 + \delta]$.
\end{lemma}

\begin{proof}
From \eqref{x-gap-2} and H\"older's inequality in time it will suffice to prove the spatial estimates
\be{s1}
\| I(fg) \|_{r'_1} \lesssim \| I f \|_{r_1} \| \langle \nabla \rangle^\eps g \|_{r_0/(p-1)}
\end{equation}
and
\be{s2}
\| I(fg) - (If) g \|_{r'_1} \lesssim N^{-\alpha} \| I f \|_{r_1} \| \langle \nabla \rangle^\eps g \|_{r_0/(p-1)}.
\end{equation}
From \eqref{x-gap} and H\"older's inequality we have
$$ \| (If)g \|_{r'_1} \lesssim \| I f \|_{r_1} \| g \|_{r_0/(p-1)}$$
so \eqref{s1} follows from \eqref{s2}.

It remains to prove \eqref{s2}.
By applying a Littlewood-Paley  decomposition to $g$, and lowering $\eps$ if necessary, we may assume that $\hat g$ is supported in the region $\langle \xi \rangle \sim M$ for some $M \geq 1$.

Fix $M$, and suppose that $\hat f$ is supported on the region $\langle \xi \rangle \lesssim M$. If $M \ll N$ then the left-hand side vanishes (since $I$ is then the identity on both $f$ and $fg$), so we may assume $M \gtrsim N$.  Then by \eqref{x-gap} and H\"older's inequality (discarding all the $I$s) we have
$$ N^\alpha \| I(fg) - (If) g \|_{r'_1}
\lesssim N^\alpha 
\| f \|_{r_1} \| g \|_{r_0/(p-1)}
\lesssim N^\alpha (M/N)^{1-s} \| If \|_{r_1} M^{-\eps} \| \langle \nabla \rangle^\eps g \|_{r_0/(p-1)}$$
which is acceptable if $s$ is close to 1 and $\alpha$ is less than $1-s$.

It remains to consider the case when $\hat f$ is supported on the region $\langle \xi \rangle \gg M$.  By dyadic decomposition we may assume that $\hat f$ is supported on the region $\langle \xi \rangle \sim 2^k M$ for some $k \gg 1$, as long as we get some exponential decay in $k$ in our estimate.

Fix $k$.  We compute the Fourier transform of $I(fg) - (If)g$:
$$ \widehat{I(fg)-(If)g}(\xi) = 
\int_{\xi = \xi_1 + \xi_2} (m(\xi_1+\xi_2) - m(\xi_1))  \hat f(\xi_1) \hat g(\xi_2).$$
From our Fourier support assumptions we may assume $\langle \xi_1 \rangle \sim 2^k M$ and $\langle \xi_2 \rangle \sim M$.  We may assume that $2^k M \gtrsim N$ since the integrand vanishes otherwise.  From the mean-value theorem we observe that
$$ m(\xi_1+\xi_2) - m(\xi_1) = O(2^{-k} m(2^k M)).$$
From this, combined with similar bounds on derivatives of $ m(\xi_1+\xi_2) - m(\xi_1)$, we obtain by standard paraproduct estimates (see e.g. \cite{christ:lectures})
$$ \| I(fg)-(If)g \|_{r'_1} \lesssim 2^{-k} m(2^k M) \| f \|_{r_1}
\| g \|_{r_0/(p-1)} \lesssim 2^{-k} \| I f \|_{r_1} M^{-\eps} 
\| \langle \nabla \rangle^\eps g \|_{r_0/(p-1)}.$$
Since $2^k M \gtrsim N$, the claim \eqref{s2} follows for $\alpha$ sufficiently small (note that we have an exponential decay in $k$ so we can safely sum in $k$).
\end{proof}

Since
$$ I \nabla F(u) = I(\nabla u \cdot F'(u))$$
we see that \eqref{chain} follows from \eqref{a}, \eqref{b} and the first part of Lemma \ref{comm}.

Now we prove \eqref{chain-2}.
By the chain rule we have
\begin{align}
\nabla (I F(u) - F(Iu)) &= I (\nabla u \cdot F'(u)) - (I \nabla u) \cdot F'(Iu) \nonumber \\
&= I (\nabla u \cdot F'(u)) - (I \nabla u) \cdot F'(u) \label{c1}\\
&\quad + (I \nabla u) \cdot (F'(u) - F'(Iu)). \label{c4}
\end{align}

The contribution of \eqref{c1} is acceptable from \eqref{a}, \eqref{b}, and the second part of Lemma \ref{comm}, so we turn to \eqref{c4}.  From Lemma \ref{strich} and \eqref{small-x} we have
$$ \| I \nabla u \|_{q_1, r_1} \lesssim 1$$
so by H\"older and \eqref{x-gap}, \eqref{x-gap-2} it suffices to show that
$$ \| F'(u) - F'(Iu) \|_{q_0/(p-1), r_0/(p-1)} \lesssim N^{-\alpha}.$$
From \eqref{small-x} and Lemma \ref{control} we have
$$ \| \langle \nabla \rangle^\beta u \|_{q_0, r_0} \lesssim 1$$
which in particular implies that
$$ \| u - Iu \|_{q_0, r_0} \lesssim N^{-\beta}.$$
The claim then follows (if $\alpha$ is sufficiently small) from \eqref{holder} and H\"older's inequality.
\end{proof}

\section{Proof of Theorem \ref{main-1}}\label{main-sec}
 
We now prove Theorem \ref{main-1}.  As in the earlier paper \cite{ckstt:7} in this series, we break the argument up into a standard series of steps.

{\bf Step 0.  Preliminaries; introduction of the modified energy.}

It suffices to show the polynomial growth bound \eqref{poly}, since the global well-posedness then follows from the local well-posedness theory in \cite{cwI}.  By another application of the local well-posedness and regularity theory and standard limiting arguments, it suffices to prove \eqref{poly} for global smooth solutions $u$ which are rapidly decreasing in space. 

Fix $u$, $s$, $T$.  We shall allow the implicit constant in $A \lesssim B$ to depend on the quantity $\|u_0\|_{H^s}$.
By time reversal symmetry we may take $T > 0$.  We will in fact assume $T \gtrsim 1$ since the case $T \ll 1$ follows from the local well-posedness theory.

Let $N \gg 1$, $\lambda \gg 1$ be quantities depending on $T$ to be chosen later.  Define the rescaled solution $u_\lambda$ by
$$ u_\lambda(x,t) := \lambda^{-2/(p-1)} u(x/\lambda, t/\lambda^2).$$
Note that $u_\lambda$ also solves the Cauchy problem \eqref{nls} but with Cauchy data
$$ (u_0)_\lambda(x) := \lambda^{-2/(p-1)} u(x/\lambda).$$

Let $I = I_N$ be as in Section \ref{notation-sec}.

We claim that
\be{hamil-init}
H(I((u_0)_\lambda)) \lesssim N^{2(1-s)} \lambda^{2(s_c-s)}.
\end{equation}
For the kinetic energy component $\| I (u_0)_\lambda \|_{\dot H^1}^2$ of the Hamiltonian, this follows from the computation
$$ \| I((u_0)_\lambda) \|_{\dot H^1}
\lesssim N^{1-s} \| (u_0)_\lambda \|_{\dot H^s} = N^{1-s}
\lambda^{s_c - s} \| u_0 \|_{\dot H^s}.$$
For the potential energy component $\| I(u_0)_\lambda \|_{p+1}^{p+1}$, we use Sobolev embedding to estimate 
$$ \| I((u_0)_\lambda) \|_{p+1} \lesssim \| I(u_0)_\lambda \|_{\dot H^{s'}}
\lesssim \| (u_0)_\lambda \|_{\dot H^{s'}}
= \lambda^{s_c-s'} \| u_0 \|_{\dot H^{s'}}
\lesssim \lambda^{s_c-s'} \| u_0 \|_{H^{s}}$$
where the exponent $s_c < s' < 1$ is determined from the Sobolev embedding theorem as 
$$ s' := \frac{n}{2} - \frac{n}{p+1}.$$
If $s$ is sufficiently close to 1, then we have
$$ (s_c-s')(p+1) \leq 2(s_c - s)$$
and so this bound is acceptable.

We now choose $N$ and $\lambda$ so that
\be{cont-1}
N^{2(1-s)} \lambda^{2(s_c-s)} \ll 1.
\end{equation}
In other words, we choose the parameters so that $I(u_0)_\lambda$ has small Hamiltonian.

In Steps 2-4 we shall prove the following almost conservation law on $H(Iu)$:

\begin{lemma}\label{almost-conserved}  Suppose that $u$ is smooth rapidly decreasing solution to \eqref{nls}, and that $H(Iu(t_0)) \lesssim 1$ for some time $t_0 \geq 0$. Then we have
\be{increment}
H(Iu(t)) = H(Iu(t_0)) + O(N^{-\alpha})
\end{equation}
for all $t_0 \leq t \leq t_0+\delta$, where $\alpha = \alpha(n,p) > 0$is a quantity depending only on $n$ and $p$, and $\delta > 0$ depends only on $n$, $p$, and the bound on $H(Iu(t_0))$.
\end{lemma}

{\bf Step 1. Deduction of \eqref{poly} from Lemma \ref{almost-conserved}.}

Let us assume Lemma \ref{almost-conserved} for the moment and deduce \eqref{poly}.

We apply the lemma to $u_\lambda$. Iterating the lemma about $N^{\alpha}$ times we obtain
$$ H(Iu_\lambda(t)) \lesssim 1 \hbox{ for all } 0 \leq t \ll N^\alpha.$$

Since we are in the defocusing case we thus have
$$ \| Iu_\lambda(t) \|_{\dot H^1} \lesssim 1 \hbox{ for all } 0 \leq t \ll N^\alpha.$$
On the other hand, by scaling the $L^2$ conservation law we have
$$ \| u_\lambda(t) \|_2 = \| u_\lambda(0) \|_2 = \lambda^{-s_c} \| u(0)\|_2 \lesssim \lambda^{-s_c}.$$
Combining these two estimates using Plancherel we obtain
$$ \| u_\lambda(t) \|_{\dot H^s} \lesssim \lambda^{-(1-s)s_c} \hbox{ for all } 0 \leq t \ll N^\alpha.$$
Undoing the scaling, we obtain
$$ \| u(t) \|_{\dot H^s} \lesssim \lambda^{s - s_c} \lambda^{-(1-s)s_c} \hbox{ for all } 0 \leq t \ll N^\alpha/\lambda^2.$$
Combining this with the $L^2$ conservation law we obtain
$$ \| u(t) \|_{H^s} \lesssim 1 + \lambda^{s - s_c} \lambda^{-(1-s)s_c} \hbox{ for all } 0 \leq t \ll N^\alpha/\lambda^2.$$
If $s$ is sufficiently close to 1, we can choose $N, \lambda \gg 1$ obeying \eqref{cont-1} such that
$$ N^\alpha/\lambda^2 \gg T.$$
The claim \eqref{poly} then follows by unraveling the exponents.

It only remains to prove Lemma \ref{almost-conserved}.  This will be done in the next three steps.

{\bf Step 2.  Control $u$ at time $t_0$.}

Let $u$ be as in Lemma \ref{almost-conserved}. The hypothesis $H(Iu(t_0)) \lesssim 1$ immediately implies that
\be{init}
\| Iu(t_0) \|_{\dot H^1} \lesssim 1.
\end{equation}

{\bf Step 3.  Control $u$ on the time interval $[t_0, t_0+\delta]$.}

We will make the a priori assumption that
$$ \sup_{t_0 \leq t \leq t_0 + \delta} H(Iu(t)) \leq C$$
for some large constant $C$; this assumption can then be removed by the usual limiting arguments.  In particular we have
\be{lp1}
\| Iu \|_{\infty, p+1} \lesssim 1
\end{equation}
since we are in the defocusing case.  Here and in the rest of this section, we adopt the convention that all spacetime norms are over the region $\R^n \times [t_0,t_0+\delta]$.

The next step is prove the spacetime estimate
\be{X-bound}
\| Iu \|_X \lesssim 1
\end{equation}
where $X$ is the space defined in Section \ref{strichartz-sec}.  By another continuity argument we may
make the \emph{a priori} assumption that $\| Iu \|_X < C$ for some large constant $C$.  We compute
\bas
\| Iu \|_X &\lesssim \| Iu(t_0)\|_{\dot H^1} + \| I \nabla (iu_t + \Delta u) \|_{q'_1, r'_1} \\
&\lesssim \| Iu(t_0)\|_{\dot H^1} + \| I \nabla F(u) \|_{q'_1, r'_1}.
\end{align*}
But then \eqref{X-bound} follows from \eqref{init} and Lemma \ref{chain-rule}, if $\delta$ is chosen sufficiently small.

From \eqref{X-bound} and \eqref{lp1} we see that \eqref{small-x} holds, so that Lemma \ref{chain-rule} is now available.

{\bf Step 4.  Control the increment of the modified energy.}

It remains to deduce \eqref{increment} from \eqref{X-bound} and \eqref{lp1}.  By the fundamental theorem of Calculus it suffices to show
$$ |\int_0^t \partial_t H(Iu)\ dt| \lesssim N^{-\alpha}.$$
From an integration by parts we have
$$ \partial_t H(Iu) = \Re \int \overline{Iu_t} (-\Delta Iu + F(Iu))\ dx.$$
Since $Iu_t = i\Delta Iu - i I F(u)$, we have
$$ \Re \int \overline{Iu_t} (-\Delta Iu + IF(u))\ dx = 0$$
and so
$$ \partial_t H(Iu) = \Re \int \overline{Iu_t} (F(Iu) - IF(u))\ dx.$$
Expanding $Iu_t$ and integrating by parts, it thus suffices to prove the estimates 
\be{nab}
|\int_0^t \overline{I\nabla u} \cdot \nabla (IF(u) - F(Iu))|
\lesssim N^{-\alpha} 
\end{equation}
and
\be{nonab}
|\int_0^t \overline{IF(u)} (IF(u) - F(Iu))| \lesssim N^{-\alpha}.
\end{equation}

We first prove \eqref{nab}. From \eqref{X-bound} and Lemma \ref{strich} we have
$$ \| I \nabla u \|_{q_1, r_1} \lesssim 1$$
and the claim then follows from H\"older's inequality and Lemma \ref{chain-rule}.

Now we prove \eqref{nonab}. By Cauchy-Schwarz we may estimate the left-hand side as
$$ \| IF(u) \|_{2,2} \| IF(u) - F(Iu) \|_{2,2}.$$
If $s$ is sufficiently close to 1 we have
$$ \| IF(u) \|_{2,2} \lesssim \| F(u) \|_{2,2} \lesssim \| u \|_{2p, 2p}^p \lesssim \| \langle \nabla \rangle^\gamma Iu \|_{2p, 2p}^p \lesssim 1$$
where the last inequality follows from \eqref{X-bound}, \eqref{lp1}, and Lemma \ref{control}.  Thus it suffices to show that
$$ \| IF(u) - F(Iu) \|_{2,2} \lesssim N^{-\alpha}.$$
Split $u = u_{low} + u_{high}$, where $u_{low}$ and $u_{high}$ are smoothly cut off projections of $u$ to the Fourier regions $|\xi| < N/2$ and $|\xi| > N/4$ respectively.  From \eqref{triangle} we have 
\bas 
F(u) &= F(u_{low}) + O(|u_{high}| |u_{low}|^{p-1}) + O(|u_{high}|^p) \\
F(Iu) &= F(u_{low}) + O(|Iu_{high}| |u_{low}|^{p-1}) + O(|Iu_{high}|^p). 
\end{align*}
Also, from Lemma \ref{control} and \eqref{X-bound} we have
\be{reg}
\| \langle \nabla \rangle^\gamma u_{low} \|_{2p, 2p} \lesssim 1; \quad \| u_{high} \|_{2p, 2p} \lesssim N^{-\gamma}.
\end{equation}
Thus the error terms $O(|u_{high}| |u_{low}|^{p-1})$, $O(|u_{high}|^p)$, $O(|Iu_{high}| |u_{low}|^{p-1})$, $O(|Iu_{high}|^p)$ are acceptable from H\"older (if $\alpha$ is sufficiently small), and so it suffices to show that
$$ \| (1-I) F(u_{low}) \|_{2,2} \lesssim N^{-\alpha}.$$
To show this we exploit the additional regularity of $u_{low}$ in \eqref{reg} and compute
\bas
\| (1-I) F(u_{low}) \|_{2,2} &\lesssim N^{-1} \| \nabla F(u_{low}) \|_{2,2} \\
&\lesssim N^{-1} \| |\nabla u_{low}| |u_{low}|^{p-1} \|_{2,2} \\
&\lesssim N^{-1} \| \nabla u_{low} \|_{2p, 2p} \| u_{low} \|_{2p, 2p}^{p-1} \\
&\lesssim N^{-1} N^{1-\gamma} \| \langle \nabla \rangle^\gamma u_{low} \|_{2p, 2p} \| u_{low} \|_{2p, 2p} \\
&\lesssim N^{-\gamma}
\end{align*}
as desired.

This concludes the proof of \eqref{increment}.  The proof of Theorem \ref{main-1} is thus complete.

\section{Proof of Theorem \ref{main-2}}\label{2-sec}

We now begin the proof of Theorem \ref{main-2}.  The idea is to modify the previous argument to use the Lyapunov functional $L$ instead of the Hamiltonian $H$.

Let $\sigma$, $u_0$, $u$, $p$ be as in Theorem \ref{main-2}.  Again, we use time reversal symmetry to restrict ourselves to the case $t>0$.

From the global well-posedness theory we know that $u(t)$ is in $H^s$ globally in time; by limiting arguments we may also assume a priori that $u(t)$ is smooth and rapidly decreasing.  From $L^2$ norm conservation we have
$$ \| u(t) \|_2 = \|u_0 \|_2$$
for all $t$.
In particular, since $u_0$ is close in $H^s$ to a ground state, we have 
\be{2-bound}
\|u(t) \|_2 = \|u_0\|_2 \lesssim 1.
\end{equation}

Let $N$ be a large number depending on $\sigma$ to be chosen later.  We consider the quantity $L(Iu(t))$.

We begin by computing this quantity at time 0.  From hypothesis we have
$$ \dist_{H^s}(u_0, \Sigma) \lesssim \sigma.$$
Since $I$ maps $H^s$ to $H^1$ with operator norm $O(N^{1-s})$, we thus have
$$ \dist_{H^1}(Iu_0, I\Sigma) \lesssim N^{1-s} \sigma.$$
Also, since the ground states in $\Sigma$ are uniformly smooth, we have
$$ \dist_{H^1}(I\Sigma, \Sigma) \lesssim N^{-100}.$$
If we thus choose $N := \sigma^{-a}$ for some $0 < a < 1/(1-s)$ we thus have that
$$ \dist_{H^1}(Iu_0, \Sigma) \ll 1.$$
From \eqref{weinstein} we thus have
$$ L(Iu(0)) \lesssim 1.$$

We now prove the following analogue of Lemma \ref{almost-conserved}:

\begin{lemma}\label{almost-conserved-2}  Suppose that $u$ is smooth rapidly decreasing solution to \eqref{nls}, and that $L(Iu(t_0)) \lesssim 1$ for some time $t_0 \geq 0$. Then we have
\be{increment-2}
L(Iu(t)) = L(Iu(t_0)) + O(N^{-\alpha})
\end{equation}
for all $t_0 \leq t \leq t_0+\delta$, where $\alpha = \alpha(n,p) > 0$ is a quantity depending only on $n$ and $p$ and $\delta > 0$ depends on $n$, $p$, and the bound on $L(Iu(t_0))$.
\end{lemma}

\begin{proof}
By continuity arguments we may make the a priori assumption that $L(Iu(t)) \lesssim 1$ for all $t_0 \leq t \leq t_0 + \delta$.
From \eqref{2-bound} we have
$$ \| Iu(t) \|_2 \lesssim 1.$$
By the Gagliardo-Nirenberg inequality and the fact that $1 < p < 1+4/n$ (which comes from the hypothesis $s_c < 0$) we then have
$$ \| Iu(t) \|_{p+1}^{p+1} \lesssim \| Iu(t) \|_{\dot H^1}^{2-\theta}$$
for some $\theta > 0$ (in fact $\theta = 2 - \frac{n(p-1)}{2}$).  Since $L(Iu(t))$ is a combination of $\| Iu(t) \|_2^2$, $\| Iu(t)\|_{p+1}^{p+1}$ and $\| Iu(t) \|_{\dot H^1}^2$, we thus have that
\be{i-bound}
\| Iu(t)\|_{p+1}^{p+1}, \| Iu(t) \|_{\dot H^1}^2 \lesssim 1
\end{equation}
for all $t_0 \leq t \leq t_0 + \delta$.

We can now run the argument in the proof of Lemma \ref{almost-conserved} to prove that
$$ H(Iu(t)) = H(Iu(t_0)) + O(N^{-\alpha});$$
admittedly we are in the focusing case rather than the defocusing case, but an inspection of the argument shows that this does not matter thanks to the bounds \eqref{i-bound}.

Since $L(Iu(t))$ is a combination of $H(Iu(t))$ and $\| Iu(t) \|_2^2$, it thus suffices to show that
$$ \| Iu(t) \|_2^2 = \| Iu(t_0) \|_2^2 + O(N^{-\alpha}).$$
From $L^2$ norm conservation it suffices to show that
$$ \| Iu(t) \|_2^2 = \| u(t)\|_2^2 + O(N^{-\alpha})$$
for all $t_0 \leq t \leq t_0 + \delta$.  Since $\|u(t)\|_2$ is bounded, this is equivalent to
$$ \| Iu(t) \|_2 = \| u(t) \|_2 + O(N^{-\alpha}).$$
But from \eqref{i-bound} we have
$$ \| u(t) - Iu(t)\|_2 \lesssim N^{-1} \| Iu(t) \|_{\dot H^1} \lesssim N^{-1}$$
and the claim follows from the triangle inequality.
\end{proof}

By iterating \eqref{increment-2} we see that $L(Iu(t))$ is bounded for all $0 \leq t \ll N^\alpha$.  By applying Gagliardo-Nirenberg as in the proof of the above Lemma, we thus see that
$$ \| Iu(t) \|_{\dot H^1} \lesssim 1 \hbox{ for all } 0 \leq t \ll N^\alpha.$$
Since $Iu(t)$ was already bounded in $L^2$, we thus have
$$ \| Iu(t) \|_{H^1} \lesssim 1 \hbox{ for all } 0 \leq t \ll N^\alpha.$$
If $Iu(t)$ is bounded in $H^1$, then $u(t)$ is bounded in $H^s$.  Thus $u(t)$ stays in a bounded ball in $H^s$ for time $t \ll N^\alpha$.  Since $N = \sigma^{-a}$, the claim follows.


\begin{thebibliography}{10}

\bibitem{cwI}
T. Cazenave, F.B. Weissler, \emph{Critical nonlinear Schr\"odinger
Equation}, Non. Anal. TMA, \textbf{14} (1990), 807--836.

\bibitem{bourg.2d}
J. Bourgain, \emph{Refinements of Strichartz Inequality and Applications
to 2d-NLS With Critical Nonlinearity}, Inter. Math. Res. Not.,  (1998),
p. 253--284.

\bibitem{borg:scatter}
J. Bourgain, \emph{Scattering in the energy space and below for 3D NLS}, J. Anal. Math. \textbf{75} (1998), 267-297. 

\bibitem{borg:book}
J. Bourgain, \emph{New global well-posedness results for non-linear Schr\"odinger equations}, AMS Publications, 1999.

\bibitem{christ:lectures}
M. Christ,
Lectures on singular integral operators. 
CBMS Regional Conference Series in Mathematics, 77. 
Published for the Conference Board of the Mathematical Sciences, Washington, DC; by the American Mathematical Society, Providence, RI, 1990. 

\bibitem{cct}
M. Christ, J. Colliander, T. Tao, \emph{Asymptotics, modulation, and low-regularity ill-posedness of canonical defocusing equations}, preprint.

\bibitem{christ:weinstein}
M. Christ, M. Weinstein, \emph{Dispersion of small amplitude solutions of the
 generalized Korteweg-de Vries equation}, J. Funct. Anal. \textbf{100} (1991), 
 87--109.

\bibitem{coff}
C.V. Coffman, \emph{Uniqueness of the ground state solution for $\Delta u - u + u^3 = 0$ and a variational characterization of other solutions}, Arch. Rat. Mech. Anal. \textbf{46} (1972), 81--95.

\bibitem{coifmanm1} Coifman, R. R and Meyer, Y. 
\emph{On commutators of singular integrals and bilinear singular integrals}
, Trans. AMS 212, pp. 315--331 [1975], 

\bibitem{coifmanm5} Coifman, R. R and Meyer, Y. 
\emph{Non-linear harmonic analysis, operator theory and P.D.E},
Beijing Lectures in Analysis, Annals of Math. Studies 112,
3-46 [1986]

\bibitem{coifmanm6} Coifman, R. R and Meyer, Y. 
\emph{Ondelettes et op\'erateurs III, Op\'erateurs multilin\'eaires},
Actualit\'es Math\'ematiques, Hermann, Paris 1991

\bibitem{ckstt:1}
J. Colliander, M. Keel, G. Staffilani, H. Takaoka, T. Tao, \emph{Almost conservation laws and global rough solutions to a nonlinear Schrodinger equation}, to appear, Math. Res. Letters. 

\bibitem{ckstt:2}
J. Colliander, M. Keel, G. Staffilani, H. Takaoka, T. Tao, \emph{Sharp global well-posedness for periodic and non-periodic KdV and mKdV on $\R$ and $\T$}, to appear, J. Amer. Math. Soc. 

\bibitem{ckstt:3}
J. Colliander, M. Keel, G. Staffilani, H. Takaoka, T. Tao, \emph{Multilinear estimates for periodic KdV equations, and applications}, to appear, J. Funct. Anal. 

\bibitem{ckstt:4}
J. Colliander, M. Keel, G. Staffilani, H. Takaoka, T. Tao, \emph{Further global well-posedness results for the 2D cubic NLS}, in preparation.

\bibitem{ckstt:5}
J. Colliander, M. Keel, G. Staffilani, H. Takaoka, T. Tao, \emph{Global well-posedness for Schr\"odinger equations with derivative}, SIAM J. Math. Anal. \textbf{33} (2001), 649--669 (electronic).

\bibitem{ckstt:6}
J. Colliander, M. Keel, G. Staffilani, H. Takaoka, T. Tao, \emph{A refined global well-posedness result for the Schrodinger equations with 
 derivative}, to appear, Siam J. Math. Anal. 2002.

\bibitem{ckstt:7}
J. Colliander, M. Keel, G. Staffilani, H. Takaoka, T. Tao, \emph{Polynomial growth and orbital instability bounds for the 1D cubic NLS below the energy norm}, in preparation.

\bibitem{ckstt:scatter}
J. Colliander, M. Keel, G. Staffilani, H. Takaoka, T. Tao, \emph{Scattering for the 3D cubic NLS below the energy norm}, in preparation.

\bibitem{feff:max}
C. Fefferman and E.~M. Stein, \emph{Some maximal inequalities}, Amer.
J. Math. \textbf{93} (1971): 107--115.

\bibitem{gv:scatter}
J. Ginibre, G. Velo, \emph{Scattering theory in the energy space for a class of nonlinear Schr\"odinger equations}, J. Math. Pure. Appl. \textbf{64} (1985), 363--401.

\bibitem{tao:keel}
M. Keel, T. Tao, \emph{Endpoint Strichartz Estimates}, Amer. Math. J. 120 (1998), 955--980.

\bibitem{keel:wavemap}
M. Keel, T. Tao, \emph{Local and global well-posedness
of wave maps on $\R^{1+1}$ for rough data}, IMRN \textbf{21} (1998),
1117--1156.

\bibitem{keel:mkg}
M. Keel, T. Tao, \emph{Global well-posedness
of the Maxwell-Klein-Gordon equation below the energy norm}, to appear.

\bibitem{kpv:nlw}
C. Kenig, G. Ponce, L. Vega, \emph{Global well-posedness for semi-linear wave equations}, preprint. 

\bibitem{kwong}
M.K. Kwong, \emph{Uniqueness of positive solutions of $\Delta u - u + u^p = 0$ in $\R^n$}, Arch. Rat. Mech. Anal. \textbf{105} (1989), 243--266.

\bibitem{mc}
K. McLeod, J. Serrin, \emph{Nonlinear Schr\"odinger equation.  Uniqueness of positive solutions of $\Delta u + f(u) = 0$ in $\R^n$}, Arch. Rat. Mech. Anal. \textbf{99} (1987), 115-145.

\bibitem{merle}
F. Merle, \emph{Determination of blow-up solutions with minimal mass for nonlinear Schrodinger equation with critical power}, Duke Math. J. 69 (1993), 427-453. 

\bibitem{nak:scatter}
K. Nakanishi, \emph{Energy scattering for non-linear Klein-Gordon and Schrodinger equations in spatial dimensions 1 and 2}, JFA \textbf{169} (1999), 201--225. 

\bibitem{stein:small}
E.~M. Stein, \emph{Singular Integrals and Differentiability
Properties of Functions}, Princeton University Press, 1970.

\bibitem{sulem}
C. Sulem, P. Sulem, The nonlinear Schrodinger equation: Self-Focusing and Wave Collapse, Applied Mathematical Sciences 139, Springer-Verlag, New York.

\bibitem{tak:dnls}
H. Takaoka, \emph{Global well-posedness for Schrodinger equations with derivative in a nonlinear term and data in low-order Sobolev spaces}, Electron. J. Diff. Eqns. \textbf{42} (2001), 1--23. 

\bibitem{tvv:bilinear}
T. Tao, A. Vargas, L. Vega, \emph{A bilinear approach to the
restriction and Kakeya conjectures},
 J. Amer. Math. Soc. \textbf{11} (1998), pp. 967--1000.

\bibitem{wein:modulate}
M. Weinstein, \emph{Modulational stability of ground states of nonlinear Schrodinger equations}, SIAM J. Math. Anal. \textbf{16} (1985), 472-491. 

\bibitem{wein}
M. Weinstein, \emph{Lyapunov stability of ground states of nonlinear dispersive equations}, CPAM \textbf{39} (1986), 51-68. 


\end{thebibliography}
\end{document}